\documentclass[11pt]{article}
\usepackage[english]{babel}
\usepackage{bm}
\usepackage[letterpaper,top=2cm,bottom=2cm,left=3cm,right=3cm,marginparwidth=1.75cm]{geometry}

\usepackage[numbers,sort&compress]{natbib}

\usepackage{amsmath}
\usepackage{url}
\allowdisplaybreaks
\usepackage{amssymb}
\usepackage{amsthm}
\usepackage{graphicx}
\usepackage[colorlinks=true, allcolors=blue]{hyperref}

\newtheorem{proposition}{Proposition}[section]
\newtheorem{theorem}{Theorem}[section]
\newtheorem{corollary}[theorem]{Corollary}
\newtheorem{lemma}[proposition]{Lemma}
\newtheorem{definition}{Definition}[section]
\newtheorem{remark}{Remark}[section]

\title{Quasi-periodic  Dynamics for Multi-dimensional Quasi-linear 
Schr\"{o}dinger  Equations via Resonant Mode Control}
\author{Zuhong You \thanks{The work was supported  by National Natural Science Foundation of China ( Grant NO. 12371189).} $^1$ and Xiaoping Yuan \footnotemark[1] $^1$  \\
		$^{1}$School of Mathematical Sciences, Fudan University, Shanghai 200433, P. R. China   }
\date{}

\begin{document}

\maketitle
\begin{abstract}

This paper focuses on the problem of quasi-periodic solutions for multi-dimensional quasi-linear Schr\"odinger equation. To address the challenge of unbounded perturbations caused by quasi-linear terms in the equation, we define the resonant mode set $\mathcal{K}$ to control nonlinear resonant effects. Combining KAM (Kolmogorov-Arnold-Moser) ( or Nash-Moser ) theory and Fourier analysis methods, we prove
that there are plenty of quasi-periodic solutions of the  equation. We also present the Fourier expansion form of the solutions and the estimation of frequency shifts. 

\end{abstract}
\section{Introduction}

The nonlinear Schr\"{o}dinger equation (NLSE) stands as a cornerstone in mathematical physics, describing phenomena ranging from quantum mechanics (wave functions of particles in nonlinear potentials) to optics (propagation of ultrashort pulses in Kerr media ) and plasma physics (dynamics of ion-acoustic waves). See \cite{Bia1976nonlinearwave, Chakrabort2021, Chancelor2017}, for example. When nonlinearities involve derivatives of the same order as the linear operator (i.e., \textit{quasi-linear} terms), the equation is called the quasi-linear 
Schr\"{o}dinger equation which 
captures more complex physical behaviors—for instance, dispersion-managed nonlinear optics (where nonlinear and dispersive effects compete at the same order ) or quantum systems with strong nonlinear couplings (e.g., Bose–Einstein condensates under high-intensity interactions). The quasi-linear 
Schr\"{o}dinger equation models the self-channeling
of a high-power ultra short laser in matter. It also appears in plasma physics and fluid mechanics, in the theory of Heisenberg ferromagnets and magnons, in dissipative quantum mechanics and in condensed
matter theory, and so on. 
See \cite{Borovskii1993, Bra1993, Kuri81, Laed83, Bass90, Kos90, Has80, Mak84}, for example.

Here we address multi-dimensional quasi-linear Schr\"odinger equations of the form:
\begin{equation}
    iu_{t}-\Delta u+\rho u+\Delta (\frac{\partial H}{\partial \bar{u}}(u,\bar{u}))=0,
    \label{yuanshifangcheng}
\end{equation}
with periodic boundary conditions ($x \in \mathbb{T}^d$). Here, $H(u, \bar{u}) = |u|^4 + \sum_{\substack{l+l' \geq 5 \\ l,l' \geq 0}} \alpha_{l,l'} u^l \bar{u}^{l'}$ is analytic near $0$, $\alpha_{l,l'} \in \mathbb{R}$, and $\alpha_{l,l'} = \alpha_{l',l}$ (ensuring symmetry). Our goal is to establish the existence of quasi-periodic solutions for \eqref{yuanshifangcheng}—solutions that exhibit ordered yet non-repeating temporal/spatial behavior, with profound implications for understanding  stable, non-chaotic dynamics in nonlinear systems (e.g., persistent coherent structures in quantum fluids). These solutions describe systems with multiple incommensurate frequencies (e.g., a quantum particle in a nonlinear potential with multiple
stable energy levels, or a nonlinear optical pulse with superposed frequency components that never fully repeat ). Their existence implies long-term stability-a critical property for applications like secure optical communications or robust quantum device design.


The construction of quasi-periodic solutions for PDEs has   emerged two  approaches: the classical KAM technique (see, e.g., \cite{Kuk87, Way90, Kuk93, Pos96, LY00, eliasson2010kam, procesi2015kam, eliasson2016kam, yuan2021kam}) and the Craig-Wayne-Bourgain (CWB) method (see, e.g., \cite{craig1993newton, bourgain1995construction, bourgain1998quasi, bourgain2005green, berti2010sobolev, BP11duke, BBnlwquasi, berti2013quasi, wang2016energy, wang2016quasi, BB20}).

When addressing unbounded perturbations, where nonlinearities involve derivative terms, KAM theory encounters significant challenges. The first breakthroughs in handling such unbounded perturbations were pioneered by Kuksin \cite{Kuk98kdv, Kuk00} who established the Kuksin Lemma  and developed applicable KAM theorems to analyze the persistence of finite-gap solutions for the  KdV equation.
Subsequently, Liu and Yuan \cite{LY10spec, LY11dnls} extended the Kuksin Lemma to its limiting case and established KAM theorems for quantum Duffing oscillators, derivative nonlinear Schrödinger equations (DNLS), and Benjamin-Ono equations. 
Later, Berti, Biasco and Procesi\cite{BBP13dnlw}  established the KAM theorems for derivative nonlinear wave equations (DNLW) by exploiting the quasi-T\"{o}plitz property of the perturbation.  Note that all these developments  rely on classical KAM techniques. Also see \cite{Sun2023} for a related problem. 
Recently, by the Craig-Wayne-Bourgain (CWB) method,  \cite{you2025periodic} and \cite{xyy2025construction}    established the existence of periodic solutions for  multi-dimensional nonlinear Schr\"{o}dinger equations (NLS) and multi-dimensional nonlinear wave  equations (NLW) with fractional derivative perturbations, respectively.

All the aforementioned results concern  PDEs in
which the nonlinearity contains strictly less derivatives than the linear differential
operator. For quasi-linear or fully nonlinear  equations, namely
PDEs whose nonlinearities contain derivatives of the same order as the linear operator, the KAM question becomes quite delicate. 
The first breakthroughs concerned the periodic solutions for fully nonlinear PDEs, where a number of pioneering works such as \cite{rabinowitz1967periodic, rabinowitz1969periodic, iooss2005standing, rard2009small, baldi2009periodic, 
iooss2011asymmetrical, baldi2013periodic},  made seminal contributions using diverse methods. For quasi-periodic solutions in fully nonlinear PDEs,  Baldi, Berti and
Montalto \cite{BBM14, BBM16kdv}  developed a novel framework to tackle the challenges and established the KAM theorem for quasi-linear  forced
Airy equation and autonomous quasi-linear  KdV. Later, Feola and Procesi \cite{feola2015quasi}
generalized their strategy to the NLS and established the KAM theorem for one-dimensional fully nonlinear forced reversible Schr\"{o}dinger equations. 
Subsequently, the framework has been extended and applied to various problems in fluid dynamics, yielding significant results.  See, for example, \cite{BBHM18, BM20fl, BFM20fl, BFM24fl, FG24, BFMT25} and the references therein.


As for the existence of quasi-periodic solutions to quasi-linear Schr\"odinger equations of spatial multi-dimension $d>1$, there is no any result at the present time, according to our knowledge.  Inspired by \cite{yuan2003quasi},  we are  devoted to proving the existence of $b$-dimensional ($b \le d$) quasi-periodic solutions for $d$-dimensional quasi-linear Schr\"{o}dinger equations. 

The unperturbed equation
\begin{equation}
    iu_{t}-\Delta u+\rho u=0,\ \  x\in \mathbb T^d,
\end{equation}
has quasi-periodic solutions
\begin{equation}
u_{0}(x,t)=\sum\limits_{j=1}^{b}a_{j}e^{in_{j}\cdot x+i\omega_{n_{j}t}},\ \  a_{j}\in\mathbb{C},
\label{feiraojie}
\end{equation}
where $n_{j}\in \mathbb{Z}^{d}$ and $\omega_{n_{j}}=|n_{j}|^{2}+\rho$ for $1\le j\le b$. Physically, $u_0$ represents a superposition of monochromatic waves (each with wavevector $n_j$ and frequency $\omega_{n_j}$), a common model for coherent nonlinear modes in quantum or optical systems (e.g., discrete-mode solitons in photonic lattices).

We investigate the \textbf{persistence} of the quasi-periodic solution \eqref{feiraojie} under the quasi-linear perturbation $\Delta\left(\frac{\partial H}{\partial \bar{u}}(u, \bar{u})\right)$.
Precisely, we say the quasi-periodic solution \eqref{feiraojie} persists with frequency $\bm{\omega}$, if there exists a solution of \eqref{yuanshifangcheng} expressed in the form 
\begin{equation}
    u(x,t)=\sum\limits_{\substack{n\in\mathbb{Z}^{d}\\k\in\mathbb{Z}^{b}}}\hat{u}(n,k)e^{in\cdot x+ik\cdot \bm{\omega} t},
\end{equation}
satisfying 
\begin{align}
    &\hat{u}(n_{j},e_{j})=a_{j},\\
    &\left(\sum\limits_{(n,k)\ne(n_{j}, e_{j})}|\hat{u}(n,k)|^{2}|(n,k)|^{2s}e^{2\sigma|(k,n)|}\right)^{1/2}<|(a_{1},...,a_{b})|^{1+\kappa},\\
    &|\bm{\omega}-\bm{\omega}_{0}|<|(a_{1},...,a_{b})|^{2},
\end{align}
where $\kappa>0$, $\bm{\omega}_{0}=(\omega_{n_{1}},...,\omega_{n_{b}})$ and $s>b/2$.

We may assume $a_j \in \mathbb{R}$ (for $1 \leq j \leq b$) without loss of generality, via a time shift when $\omega_1, \dots, \omega_b$ are rationally independent.
 To state our results, we first define \textit{resonant modes}
\begin{definition}[Resonant Modes]
    Given $b$ linearly independent vectors $n_{1},...,n_{b}\in \mathbb{Z}^{d} (b\le d)$. Call
    \[ \mathcal{K}(n_{1},...,n_{b}):=\left\{k\in\mathbb Z^b: \, \sum_{j=1}^{b}k_{j}=1, \; \text{and}\; \sum_{j=1}^{b}k_{j}|n_{j}|^{2}=| \sum_{j=1}^{b} k_{j}n_{j} |^{2} \right\}. \]
  resonant modes corresponding to $n_{1},...,n_{b}$.
    \end{definition}
We have $|\mathcal{K}(n_{1},...,n_{b})|<\infty$ (see Section \ref{structureofR} for detail).
Moreover, it is easy to check that $e_{j}\in\mathcal{K}(n_{1},...,n_{b})$ for $1\le j\le b$.
Denote $\mathcal{K}_{1}=\{e_{1},...,e_{b}\}$ and $\mathcal{K}_{2}(n_{1},...,n_{b})=\mathcal{K}(n_{1},...,n_{b})\setminus\mathcal{K}_{1}$. 
Resonances ($k \in \mathcal{K}$) correspond to nonlinear interactions that \textit{synchronize} frequency components. Avoiding non-trivial resonances (i.e., $\mathcal{K}_2 = \emptyset$ in Theorem 1.1) ensures these interactions do not destroy the quasi-periodic structure, analogous to how "non-resonant" perturbations preserve stability in classical mechanics (a concept foundational to understanding persistent patterns in nonlinear optics  and quantum fluids).

We have the following statement.

\begin{theorem}\label{dingli1.1}
    Suppose that $\rho$ satisfies Diophantine conditions
    \begin{equation}
        |m\rho-l|\ge \frac{\gamma}{m^{\tau}}, \ \ \forall m\ne0, 
    \end{equation}
    for some $\gamma>0$ and $\tau>1$.  Select linearly independent $n_{1},...,n_{b}\in \mathbb{Z}^{d}$ ($d\ge b$) which satisfies that $$\mathcal{K}_{2}(n_{1},...,n_{b})=\emptyset.$$
    If $\varepsilon$ is sufficiently small, for all $(a_{1},...,a_{b})\in[\varepsilon,2\varepsilon]^{b}$, the quasi-periodic (periodic) solution \eqref{feiraojie} persists with frequency
    \begin{equation}
        \omega_{j}=|n_{j}|^{2}+\rho+|n_{j}|^{2}\left(4\sum_{l=1}^{b}a_{l}^{2}-2a_{j}^{2}\right)+O(|(a_{1},...,a_{b})|^{3}),\ \ j=1,...,b.\label{piaoyipinlv}
    \end{equation}
\end{theorem}

The perturbed frequency $\omega_j$'s capture how nonlinearities renormalize the energy/frequency of each mode.  This reflects the effect of particle-particle interactions on collective mode frequencies, e.g., in Bose–Einstein condensates.

Theorem \ref{dingli1.1} applies to a wide range of Fourier indices $n_{1},...,n_{b}$, from which we deduce the following corollary.
\begin{corollary}\label{tuilun1.2}
    Suppose that $\rho$ satisfies Diophantine conditions
    \begin{equation}
        |m\rho-l|\ge \frac{\gamma}{m^{\tau}}, \ \ \forall m\ne0, 
    \end{equation}
    for some $\gamma>0$ and $\tau>1$.
    \begin{itemize}
        \item[\textup{(i)}] All periodic solutions 
        \begin{equation}
            u_{0}(x,t)=ae^{in\cdot x+i\omega_{n}t}
        \end{equation}
        persist with frequency
        \begin{equation}
            \omega=|n|^{2}+\rho+2a^{2}|n|^{2}+O(|a|^{3}),
        \end{equation}
        as long as $a\in\mathbb{R}$ and $|a|$ is small enough.
        \item[\textup{(ii)}] For $d\ge 2$, select $n_{1}, n_{2}\in\mathbb{Z}^{d}$ such that $n_{1}\nparallel n_{2}$. The quasi-solution \eqref{feiraojie} persists with frequency  \eqref{piaoyipinlv}  for all $(a_{1},a_{2})\in[\varepsilon,2\varepsilon]^{2}$ as long as $\varepsilon$ is small enough.
        \item[(\textup{iii})] For $d\ge 3$, select linearly independent $n_{1},n_{2}, n_{3}$ such that the angle between $n_{i}-n_{j}$ and $n_{i}-n_{l}$ is less than $\frac{\pi}{2}$ for all $(i,j,l)\in\{\textup{the permutations of }(1,2,3)\}$. The quasi-solution \eqref{feiraojie} persists with frequency  \eqref{piaoyipinlv} for all  $(a_{1},a_{2}, a_{3})\in[\varepsilon,2\varepsilon]^{3}$ as long as $\varepsilon$ is small enough.
        \item[\textup{(iv)}] Suppose that $|n_{j}|^{2}\equiv M$, where $M>0$ is a constant, and $n_{j}\perp n_{l}$ for all $j\ne l$. The quasi-solution \eqref{feiraojie} persists with frequency  \eqref{piaoyipinlv} for all  $(a_{1},..., a_{b})\in[\varepsilon,2\varepsilon]^{b}$ as long as $\varepsilon$ is small enough.
    \end{itemize}
\end{corollary}

We now consider the case  $\mathcal{K}_{2}(n_{1},...,n_{b})\ne0$. To state our theorem, we define the following matrix.
\begin{definition}
    For $k,k'\in\mathcal{K}_{2}(n_{1},...,n_{b})$ and $k\ne k'$,
    define 
    \begin{equation}
        A_{k,k'}=2\sum\limits_{(j,l)\in J_{1}}a_{j}a_{l}+\sum\limits_{(j,l)\in J_{2}}a_{j}a_{k},
    \end{equation}
    where
\begin{align}
    &J_{1}=\{(j,l)|k+e_{i}-e_{j}=k'\},\notag\\
    &J_{2}=\{(j,l)|e_{i}+e_{j}-k=k'\}.\notag
\end{align}

For $k\in\mathcal{K}_{2}(n_{1},...,n_{b})$, define
\begin{equation}
    A_{kk}=\sum_{j=1}^{b}k_{j}|n_{j}|^{2}a_{j}^{2}.
\end{equation}
\end{definition}
\begin{remark}
    Note that $A$ is determined by $n_{1},...,n_{b}$, and $A_{kk'}\in\mathbb{Z}[a_{1},...,a_{b}]$.
\end{remark}
\begin{theorem}\label{dingli1.3}
    Suppose that $\rho$ satisfies Diophantine conditions
    \begin{equation}
        |m\rho-l|\ge \frac{\gamma}{m^{\tau}}, \ \ \forall m\ne0, 
    \end{equation}
    for some $\gamma>0$ and $\tau>1$.
    Select linearly independent $n_{1},...,n_{b}\in\mathbb{Z}^{d}$ ($d\ge b$) such that 
    \begin{itemize}
        \item[\textup{(i)}] $(n_{i}-n_{j})\not\perp(n_{i}-n_{l})$ for all distinct $i,j,l\in\{1,...,b\}$,
        \item[\textup{(ii)}] $\det A\not\equiv 0$. 
    \end{itemize}
    For $\varepsilon$ small enough, there exists a subset $I_{\varepsilon}\subset[\varepsilon,2\varepsilon]^{b}$. For all $(a_{1},...,a_{b})\in[\varepsilon,2\varepsilon]^{b}\setminus I_{\varepsilon}$, the solution \eqref{feiraojie} persists with frequency \eqref{piaoyipinlv}. Moreover, we have
    \begin{equation}
        \textup{mes } I_{\varepsilon}<\varepsilon^{b+c},
    \end{equation}
    where $c$ is a constant depending on $n_{1},...,n_{b}$.
\end{theorem}
If $H(u,\bar{u})$ takes the form
\begin{equation}
    H(u,\bar{u})=|u|^{4}+\sum_{p=3}^{\infty}\alpha_{p}|u|^{2p},
\end{equation}
where $\alpha_{p}\in \mathbb{R}$.
We have the following statement.
\begin{theorem}\label{dingli1.4}
    Theorem \ref{dingli1.1}, Corollary \ref{tuilun1.2} and Theorem \ref{dingli1.3}  hold  for all $\rho\in\mathbb{R}$.
\end{theorem}

\textbf{Notations.}
$a\lesssim b$ refers to $a\leq Cb$, for some constant which can be chosen depending on the context. $a\sim b$ refers to $b\lesssim a\lesssim b$.

$k$ is reserved to represent vectors in $\mathbb{Z}^{b}$, and $k_{j}$ refers to the $j$-th component of $k$.

\section{Lyapunov–Schmidt Reduction}
Finding a quasi-periodic solution with frequency $\omega\in\mathbb{R}^{b}$ to \eqref{yuanshifangcheng} is equivalent to finding a solution $u\in L^{2}(\mathbb{T}^{b+d})$ to
\begin{equation}
    i\partial_{\omega}u-\Delta u+\rho u +\Delta(\frac{\partial H}{\partial \bar{u}}(u,\bar{u}))=0,
    \label{hulleq}
\end{equation}
where $\partial_{\omega}:=\sum\limits_{j=1}^{b}\omega_{j}\partial_{\theta_{j}}$.
Replacing $u$ by $\varepsilon u$, we obtain
\begin{equation}
    i\partial_{\omega}u-\Delta u+\rho u+\Delta\left(\varepsilon^{-1}\frac{\partial H}{\partial \bar{u}}(\varepsilon u,\varepsilon\bar{u})\right)=0.
    \label{xiaohulleq}
\end{equation}
Notice that $\varepsilon^{-1}\frac{\partial H}{\partial \bar{u}}(\varepsilon u,\varepsilon\bar{u})=O(\varepsilon^{2})$.

Note that we have
\begin{equation}
    \textup{supp } (\frac{\partial H}{\partial \bar{u}}(u,\bar{u}))^{\land}\subset \{\sum\limits_{j=1}^{b}k_{j}(e_{j},n_{j}): k=(k_{1},...,k_{b})\in\mathbb{Z}^{b}\}:=\mathcal{W},
\end{equation}
if we have
\begin{equation}
    \textup{supp }\hat{u}\subset \mathcal{W}. 
\end{equation}
Since we aim to study the persistence of 
\begin{equation}
    u_{0}=\sum\limits_{j=1}^{b} a_{j} e^{i e_{j}\cdot\theta+in_{j}\cdot x},
\end{equation}
it is natural to find solutions to \eqref{xiaohulleq} in form
\begin{equation}
    u(\theta,x)=\sum\limits_{(k,n)\in\mathcal{W}}\hat{u}(k,n)e^{i(k,n)\cdot (\theta,x)}.
\end{equation}
Note that for $(k,n)\in\mathcal{W}$, when $k$ is fixed, $n=\sum\limits_{j=1}^{b}k_{j}n_{j}$ is fixed. We denote $n_{k}=\sum\limits_{j=1}^{b}k_{j}n_{j}$.

Let 
\begin{equation}
    H_{\mathcal{W},\sigma, s}=\left\{ u\in L^{2}(\mathbb{T}^{b+d}): \textup{supp }\hat{u}\subset \mathcal{W}, \|u\|_{\sigma,s}^{2}=C\sum\limits_{k\in\mathbb{Z}^{b}}|\hat{u}(k,n_{k})|^{2}|k|^{2s}e^{2\sigma |k|} \right\},
\end{equation}
where $\sigma>0, s>b/2$. Note that $H_{\mathcal{W},\sigma, s}$ is a Banach algebra if $C$ is chosen appropriately (see Appendix E in \cite{Berti07} for example).

Define
\begin{align}
    &T_{\omega}(u)=(i\partial_{\omega}-\Delta+\rho)u,\\
    &S(u)=\frac{\partial H}{\partial \bar{u}}( u,\bar{u}).
\end{align}
Thus, \eqref{xiaohulleq} can be written in the following form
\begin{equation}
    T_{\omega}(u)+\Delta(\varepsilon^{-1} S(\varepsilon u))=0.
    \label{fanhanform}
\end{equation}
\begin{lemma}\label{Sdexingzhi}
    Let $\sigma>0$ and $s>b/2$. Then, $S$ maps $H_{\mathcal{W},\sigma,s}$ to $H_{\mathcal{W},\sigma,s}$, and
    \begin{equation}
        \|S(u)\|_{\sigma,s}\le C_{1} \|u\|_{\sigma,s}^{3},
        \label{Sfanshu}
    \end{equation}
    for $\|u\|_{\sigma,s}$ small enough.
    Moreover, we have
    \begin{align}    \|D_{u}S(u)\|_{H_{\mathcal{W},\sigma,s}\to H_{\mathcal{W},\sigma,s}}\le C_{1} \|u\|_{\sigma,s}^{2},\label{uder}\\ \|D_{\bar{u}}S(u)\|_{H_{\mathcal{W},\sigma,s}\to H_{\mathcal{W},\sigma,s}}\le C_{1} \|u\|_{\sigma,s}^{2},\label{ubarder}
    \end{align}
    for $\|u\|_{\sigma,s}$ small enough. Here, $C_{1}$  is a constant depending on $H$, and $\|\cdot\|_{H_{\mathcal{W},\sigma,s}\to H_{\mathcal{W},\sigma,s}}$ refers to the operator norm from $H_{\mathcal{W},\sigma,s}$ to $H_{\mathcal{W},\sigma,s}$.
\end{lemma}
\begin{proof}
    Suppose $\textup{supp }\hat{u}\subset \mathcal{W}$. We have known $\textup{supp }(\frac{\partial H}{\partial \bar{u}}(u,\bar{u}))^{\land}\subset\mathcal{W}$. Moreover, we have
    \begin{align}
        \left\|\frac{\partial H}{\partial \bar{u}}(u,\bar{u})\right\|_{\sigma,s}&=\left\|2u^{2}\bar{u}+\sum\limits_{\substack{l+l'\ge5\\l\ge0,l'\ge 1}}l'\alpha_{l,l'}u^{l}\bar{u}^{l'-1}\right\|_{\sigma,s}\notag\\
        &\le 2\|u\|_{\sigma,s}^{3}+\sum\limits_{\substack{l+l'\ge5\\l\ge0,l'\ge 1}}|l'\alpha_{l,l'}|\|u\|_{\sigma,s}^{l+l'-1}\notag\\
        &\le C_{1} \|u\|_{\sigma,s}^{3},
    \end{align}
    which implies \eqref{Sfanshu}.
        On the other hand, we have
        \begin{align}
        \left\|\frac{\partial^{2} H}{\partial u\partial \bar{u}}(u,\bar{u})\right\|_{\sigma,s}&=\left\|4u\bar{u}+\sum\limits_{\substack{l+l'\ge5\\l\ge1,l'\ge 1}}ll'\alpha_{l,l'}u^{l-1}\bar{u}^{l'-1}\right\|_{\sigma,s}\notag\\
        &\le 2\|u\|_{\sigma,s}^{2}+\sum\limits_{\substack{l+l'\ge5\\l\ge1,l'\ge 1}}|ll'\alpha_{l,l'}|\|u\|_{\sigma,s}^{l+l'-2}\notag\\
        &\le C_{1} \|u\|_{\sigma,s}^{2}.\label{eq32}
    \end{align}
Since $H_{\mathcal{W},\sigma,s}$ is a Banach algebra, the estimate \eqref{eq32} implies  \eqref{uder}. Similarly, we can obtain \eqref{ubarder}.
\end{proof}
Passing \eqref{xiaohulleq} to Fourier coefficients, we have
\begin{equation}
    (-k\cdot\omega+|n|^{2}+\rho)\hat{u}(k,n)+|n|^{2}(\varepsilon^{-1}\frac{\partial H}{\partial \bar{u}}(\varepsilon u,\varepsilon\bar{u}))^{\land}(k,n)=0,\ \  (k,n)\in\mathcal{W}.
\end{equation}
Since we aim to find quasi-periodic solution with frequency near $\omega_{0}=(\omega_{n_{1}},...,\omega_{n_{b}})$, the resonant set is 
\begin{equation}
    \mathcal{R}=\{(k,n)|-k\cdot\omega_{0}+|n|^{2}+\rho=0\}\cap \mathcal{W}.
\end{equation}
Let $P$ be the projector from $H_{\mathcal{W},\sigma,s}$ to $P(H_{\mathcal{W},\sigma,s})$  defined by 
\begin{equation}                
    Pu=\sum\limits_{(k,n)\in\mathcal{W}\setminus\mathcal{R}}\hat{u}(k,n)e^{i(k,n)\cdot(\theta,x)}.
\end{equation}
Let $Q=1-P$ be defined by 
\begin{equation}
    Qu=\sum\limits_{(k,n)\in\mathcal{R}}\hat{u}(k,n)e^{i(k,n)\cdot(\theta,x)}.
\end{equation}
Then, \eqref{fanhanform} can be decomposed into 
\begin{align}
    QT_{\omega}Qu+Q\Delta(\varepsilon^{-1}S(\varepsilon Qu+\varepsilon Pu))=0,\label{bifurcation}\\
    PT_{\omega}Pu+P\Delta(\varepsilon^{-1}S(\varepsilon Qu+\varepsilon Pu))=0.\label{range}
\end{align}
The equation \eqref{bifurcation} is called the bifurcation equation and the equation \eqref{range} is called the range equation.

We adopt the following strategy to solve the equations. First, we fix $$Qu=\sum\limits_{j=1}^{b}a_{j}e^{i(e_{j},n_{j})\cdot(\theta,x)}+\sum\limits_{\substack{(k,n_{k})\notin\{(e_{j},n_{j})\}\\(k,n_{k})\in\mathcal{R}}}a_{k}e^{i(k,n_{k})\cdot(\theta,x)}$$
and $\omega$ (we also denote $a_{e_{j}}=a_{j}$). We then solve the range equation to obtain $Pu$, which is continuously differentiable with respect to $Qu$ and $\omega$. Next, we substitute $Pu$ into the bifurcation equation. By fixing $a_{j}$, we solve for $\omega$ and $a_{k}$.

\section{The Range equation}
Let 
\begin{align}
    &Qu=\sum\limits_{j=1}^{b}a_{j}e^{i(e_{j},n_{j})\cdot(\theta,x)}+\sum\limits_{\substack{(k,n_{k})\notin\{(e_{j},n_{j})\}\\(k,n_{k})\in\mathcal{R}}}a_{k}e^{i(k,n_{k})\cdot(\theta,x)}:=v,\\
    &\omega\in B(\omega_{0},\varepsilon),\label{omegafanwei}\\
    &a_{j}\in [1,2] \textup{ for } 1 \le j\le b,\label{ajfanwei}\\
    &a_{k}\in[-\varepsilon,\varepsilon]\textup{ for } (k,n_{k})\in\mathcal{R}, \textup{ where } k\ne e_{j} \textup{ for all } j.\label{akfanwei}
\end{align}

We aim to solve 
\begin{equation}
PT_{\omega}Pu+P\Delta(\varepsilon^{-1}S(\varepsilon v+\varepsilon Pu))=0.
\end{equation}

\begin{lemma}\label{baochizhengzexing}
    Suppose that $n_{1},...,n_{b}$ are linearly independent and $\rho$ satisfies Diophantine conditions
    \begin{equation}
        |m\rho-l|\ge \frac{\gamma}{m^{\tau}}, \ \ \forall m\ne0, \label{DC}
    \end{equation}
    for some $\gamma>0$ and $\tau>1$. Then, for $\omega\in B(\omega_{0},\varepsilon)$ and $\varepsilon$ sufficiently small, we have
    \begin{equation}
        \|(PT_{\omega}P)^{-1}P\Delta\|_{H_{\mathcal{W},\sigma,s}\to H_{\mathcal{W},\sigma,s}}<\frac{K(n_{1},...,n_{b},|\rho|)}{\gamma},
        \label{xianxingsuanzini}
    \end{equation}
    where $K(n_{1},...,n_{b},|\rho|)$ is a constant only depending on $n_{1},...,n_{b},|\rho|$.
\end{lemma}
\begin{proof}
    Note that $(PTP)^{-1}P\Delta$ restricted to $H_{\mathcal{W},\sigma,s}$ (modulo a Fourier transform) is
    \begin{equation}
        \textup{diag }\left(\frac{|n_{k}|^{2}}{-k\cdot\omega+|n_{k}|^{2}+\rho}:k\in\mathbb{Z}^{b}\right).
    \end{equation}
Consider the denominator $-k\cdot \omega+|n_{k}|^{2}+\rho$. Since $n_{1},...,n_{b}$ are linearly independent, we have 
\begin{equation}
    |n_{k}|^{2}=\left|\sum\limits_{j=1}^{b}k_{j}n_{j}\right|^{2}\ge c_{1} |k|^{2}.
\end{equation}
Here, $c_{1}$ is a constant depending on $n_{1},...,n_{b}$. In fact, $c_{1}=\min\limits_{\substack{|\theta|=1\\ \theta\in \mathbb{R}^{b}}} |\sum\limits_{j=1}^{b}\theta_{j}n_{j}|^{2}$. Thus, we have
\begin{align}
    |-k\cdot\omega+|n_{k}|^{2}+\rho|&\ge c_{1}|k|^{2}-|k|\max_{1\le j\le b} (|n_{j}|^{2}+|\rho|)-|\rho|\notag\\
    &\ge \frac{c_{1}}{2}|k|^{2},
\end{align}
for $|k|>\frac{\max_{1\le j\le b} (|n_{j}|^{2}+2|\rho|)}{c_{1}}:=K_{1}(n_{1},...,n_{b},|\rho|)$.

Thus, we obtain
\begin{align}
    \left|\frac{|n_{k}|^{2}}{-k\cdot\omega+|n_{k}|^{2}+\rho}\right|&\le
    \frac{\sum\limits_{1\le j\le b}|n_{j}|^{2}|k|^{2}}{\frac{c_{1}}{2}|k|^{2}}\notag\\
    &\le \frac{2\sum\limits_{1\le j\le b}|n_{j}|^{2}}{c_{1}}:=K_{2}(n_{1},...,n_{b},|\rho|),
\end{align}
for $|k|>K_{1}(n_{1},...,n_{b},\rho)$.

For $|k|\le K_{1}(n_{1},...,n_{b},|\rho|)$, we have
\begin{align}
    &-k\cdot \omega+|n_{k}|^{2}+\rho\notag\\
    =&-k\cdot\omega_{0}+|n_{k}|^{2}+\rho+O(|k|\varepsilon)\notag\\
    =&(1-\sum\limits_{j=1}^{b}k_{j})\rho+|n_{k}|^{2}-\sum\limits_{j=1}^{b}k_{j}|n_{j}|^{2}+O(|k|\varepsilon).
\end{align}
Note that $(k,n_{k})\notin \mathcal{R}$. For $1-\sum_{j=1}^{b}k_{j}=0$, we have
\begin{align}
||n_{k}|^{2}-&\sum_{j=1}^{b}k_{j}|n_{j}|^{2}|\ge 1.
\end{align}
Thus, we have
\begin{equation}
    |-k\cdot\omega+|n_{k}|^{2}+\rho|\ge 1-O(|k|\varepsilon)>\frac{1}{2},
\end{equation}
for $\varepsilon$ sufficiently small.
For $1-\sum_{j=1}^{b}k_{j}\ne0$, we have
\begin{align}
    |-k\cdot\omega+|n_{k}|^{2}+\rho|&\ge\frac{\gamma}{|1-\sum_{j=1}^{b}k_{j}|^{\tau}}-O(|k|\varepsilon)\notag\\
    &\ge \frac{\gamma}{(|k|+1)^{\tau}}-O(|k|\varepsilon)\notag\\
    &\ge \gamma K_{1}^{\tau}-O(K_{1}\varepsilon)\notag\\
    &\ge \frac{1}{2}\gamma K_{1}^{-\tau},
\end{align}
for $\varepsilon$ sufficiently small.
Thus, for $|k|\le K_{1}$, we have
\begin{equation}
    \left|\frac{|n_{k}|^{2}}{-k\cdot\omega+|n_{k}|^{2}+\rho}\right|<\frac{\sum\limits_{1\le j\le b}|n_{j}|^{2}|k|^{2}}{\frac{1}{2}\gamma K_{1}^{-\tau}}<\frac{2}{\gamma}K_{1}^{\tau+2}\sum\limits_{1\le j\le b}|n_{j}|^{2}.
\end{equation}
Denote $K(n_{1},...,n_{b}, |\rho|)=\max\left\{K_{2}(n_{1},...,n_{b},|\rho|), 2K_{1}^{\tau+2}\sum\limits_{1\le j\le b}|n_{j}|^{2}\right\}$. We complete the proof.
\end{proof}

\begin{lemma}\label{rangesol}
    Suppose  $\rho$ satisfies \eqref{DC} and $$v=\sum\limits_{j=1}^{b}a_{j}e^{i(e_{j},n_{j})\cdot(\theta,x)}+\sum\limits_{\substack{(k,n_{k})\notin\{(e_{j},n_{j})\}\\(k,n_{k})\in\mathcal{R}}}a_{k}e^{i(k,n_{k})\cdot(\theta,x)}.$$
For $\omega$, $a_{j}$, $a_{k}$ satisfying \eqref{omegafanwei}, \eqref{ajfanwei}, \eqref{akfanwei} and $\varepsilon$ sufficiently small, the equation
\begin{equation}
    PT_{\omega}Pu+P\Delta(\varepsilon^{-1}S(\varepsilon v+\varepsilon u))=0
    \label{rangeequationfanhanform}
\end{equation}
has a solution $u$, continuously differentiable with respect to  $\omega$, $a_{j}$ and $a_{k}$, such that
\begin{itemize}
    \item[\textup{(i)}] $\textup{supp } \hat{u}\subset\mathcal{W}\setminus\mathcal{R}$.
    \item[\textup{(ii)}] $\hat{u}(k,n_{k})\in\mathbb{R}$.
    \item[\textup{(iii)}] $\|u\|_{\sigma,s}\le K_{3}\varepsilon^{2}$, where $K_{3}$ is a constant depending on $n_{1},...,n_{b},\gamma$.
    \item[\textup{(iv)}] $\|\partial_{(\omega,\bm{a})}u\|_{\sigma,s}\le K_{4}\varepsilon^{2}$, where $K_{4}$ is a constant depending on $n_{1},...,n_{b},\gamma$.
\end{itemize}
\end{lemma}
\begin{proof}
    Finding a solution of \eqref{rangeequationfanhanform} is equivalent to finding a solution of
    \begin{equation}
        u=-(PT_{\omega}P)^{-1}P\Delta(\varepsilon^{-1}S(\varepsilon v+\varepsilon u)).\label{rangeequationbudongdianform}
    \end{equation}
    Denote $\mathcal{G}(u)=-(PT_{\omega}P)^{-1}P\Delta(\varepsilon^{-1}S(\varepsilon v+\varepsilon u))$.
    We next prove that $\mathcal{G}$ is a contraction map on 
    \begin{equation}
      B_{\mathbb{R}}(0, K_{3}\varepsilon^{2})=\{u\in P(H_{\mathcal{W,\sigma,s}})|\hat{u}(k,n_{k})\in\mathbb{R}, \|u\|_{\sigma,s}<K_{3}\varepsilon^{2}\},  
    \end{equation}
   where $K_{3}$ is a constant to be specified later.

    Given that $\hat{v}(k,n)\in\mathbb{R}$, it follows  that if $\hat{u}(k,n_{k})\in\mathbb{R}$ for all $k$, then
    \begin{equation}
        (\mathcal{G}(u))^{\land}(k,n)\in\mathbb{R}.
    \end{equation}
    On the other hand, since $\textup{supp }\hat{u}\subset\mathcal{W}$ and $\textup{supp }\hat{v}\subset\mathcal{W}$, we have $\textup{supp } (\varepsilon^{-1}\Delta S(\varepsilon v+\varepsilon u))\subset\mathcal{W}$. Thus, we have $\textup{supp } (\mathcal{G}(u))^{\land}\subset\mathcal{W}\setminus\mathcal{R}$.
    Let $\|u\|_{\sigma,s}<\varepsilon$. By Lemma \ref{Sdexingzhi}, we have
    \begin{align}
         \|\mathcal{G}(u)\|_{\sigma,s}\le \|(PT_{\omega}P)^{-1}P\Delta\|_{H_{\mathcal{W},\sigma,s}\to H_{\mathcal{W},\sigma,s}} \cdot \varepsilon^{-1}\cdot C_{1}\|\varepsilon v+\varepsilon u\|_{\sigma,s}^{3}\le \frac{K}{\gamma} C_{1} C_{2}\varepsilon^{2}, 
    \end{align}
       where $C_{2}$ is a constant depending on $n_{1},...,n_{b}$. Denote $K_{3}=\frac{C_{1} C_{2}K}{\gamma}$. Note that $K_{3}\varepsilon^{2}<\varepsilon$ if $\varepsilon$ is small enough. Thus, we obtain that $\mathcal{G}$ maps $B_{\mathbb{R}}(0,K_{3}\varepsilon^{2})$ to itself.

Let $u_{1}, u_{2}\in B_{\mathbb{R}}(0,K_{3}\varepsilon^{2})$. We have
\begin{align}
    &\|\mathcal{G}(u_{1})-\mathcal{G}(u_{2})\|_{\sigma,s}\notag\\
    =&\varepsilon^{-1}\|(PT_{\omega P})^{-1}P\Delta (S(\varepsilon v+\varepsilon u_{1})-S(\varepsilon v+\varepsilon u_{2}))\|_{\sigma,s}\notag\\
    \le&\varepsilon^{-1}\|(PT_{\omega P})^{-1}P\Delta\|_{H_{\mathcal{W},\sigma,s}\to H_{\mathcal{W},\sigma,s}}\notag\\
    &\cdot\left\| \int_{0}^{1}D_{u}S(\varepsilon v+\varepsilon u_{1}+\varepsilon t(u_{2}-u_{1}))\cdot\varepsilon(u_{2}-u_{1})+D_{\bar{u}}S(\varepsilon v+\varepsilon u_{1}+\varepsilon t(u_{2}-u_{1}))\cdot\varepsilon(\bar{u}_{2}-\bar{u}_{1})dt\right\|_{\sigma,s}\notag\\
    \le&\frac{K}{\gamma} C_{1}C_{3}\varepsilon^{2}\|u_{2}-u_{1}\|_{\sigma,s},
\end{align}
where $C_{3}$ is a constant depending on $n_{1},...,n_{b}$. For $\varepsilon$ sufficiently small, we have
\begin{equation}
    \|\mathcal{G}(u_{1})-\mathcal{G}(u_{2})\|_{\sigma,s}<\frac{1}{2}\|u_{2}-u_{1}\|_{\sigma,s}.
\end{equation}
Thus, there exists an unique solution $u(\omega,\bm{a})$ solving \eqref{rangeequationbudongdianform} and satisfying (i), (ii), (iii).

To analyze the regularity, let $q=\Re u$ and $p=\Im u$. Then, \eqref{rangeequationbudongdianform} becomes
\begin{align}
    q-\frac{\mathcal{G}(q+ip)+\bar{\mathcal{G}}(q+ip)}{2}=0,\label{shibu}\\
    p-\frac{\mathcal{G}(q+ip)+\bar{\mathcal{G}}(q+ip)}{2}=0.\label{xvbu}
\end{align}
Denote the left side of \eqref{shibu} and \eqref{xvbu} by $F(q,p,\omega,\bm{a})$. Note that $F$ is a map from $H_{\mathcal{W},\sigma,s}\oplus H_{\mathcal{W},\sigma,s}:=\tilde{H}_{\mathcal{W},\sigma,s}$ to itself. The norm on $\tilde{H}_{\mathcal{W},\sigma,s}$ is defined by $\|(q,p)\|_{\sigma,s}:=\|q\|_{\sigma,s}+\|p\|_{\sigma,s}$.
The equation \eqref{rangeequationbudongdianform} is equivalent to 
\begin{equation}
    F(q,p,\omega,\bm{a})=0.
\end{equation}
By \eqref{uder}, \eqref{ubarder}, \eqref{xianxingsuanzini}, we have
\begin{equation}
    \left\|\left(\frac{\partial F(q,p,\omega,\bm{a})}{\partial(q,p)}\right)^{-1}\right\|_{\tilde{H}_{\mathcal{W},\sigma,s}\to \tilde{H}_{\mathcal{W},\sigma,s}}<2,
    \label{gereenfunctionfanshuguji}
\end{equation}
for $\varepsilon$ small enough.
Note that
\begin{equation}
    \partial_{\omega_{j}} T_{\omega}=\left(\textup{diag }(-k_{j})\right)^{\lor}.
\end{equation}
We have
\begin{equation}
    \partial_{\omega_{j}} \mathcal{G}(q+ip)=-(PT_{\omega}P)^{-1}(P\partial_{\omega_{j}}TP)(PT_{\omega}P)^{-1}P\Delta(\varepsilon^{-1}S(\varepsilon v+\varepsilon u)).
\end{equation}
Following the same argument as the proof of Lemma \ref{baochizhengzexing}, we have
\begin{equation}
    \|(PT_{\omega}P)^{-1}(P\partial_{\omega_{j}}TP)\|_{H_{\mathcal{W},\sigma,s}\to H_{\mathcal{W},\sigma,s}}<\frac{K(n_{1},...,n_{b},|\rho|)}{\gamma}.
\end{equation}
Thus, we have
\begin{equation}
    \|\partial_{\omega_{j}} F(q,p,\omega,\bm{a})\|_{\sigma,s}\lesssim \left(\frac{K}{\gamma}\right)^{2}\varepsilon^{2}.
    \label{fanshuguji}
\end{equation}
By \eqref{gereenfunctionfanshuguji}, \eqref{fanshuguji} and the implicit function theorem, we obtain
\begin{align}
    \left\|\frac{\partial (q,p)}{\partial\omega_{j}}\right\|_{\sigma,s}=&\left\|\left(\frac{\partial F(q,p,\omega,\bm{a})}{\partial(q,p)}\right)^{-1} \partial_{\omega_{j}} F(q,p,\omega,\bm{a}) \right\|_{\sigma,s}\notag\\
    \le & \left\|\left(\frac{\partial F(q,p,\omega,\bm{a})}{\partial(q,p)}\right)^{-1}\right\|_{\tilde{H}_{\mathcal{W},\sigma,s}\to \tilde{H}_{\mathcal{W},\sigma,s}}\|\partial_{\omega_{j}} F(q,p,\omega,\bm{a})\|_{\sigma,s}\notag\\
    \lesssim & \left(\frac{K}{\gamma}\right)^{2} \varepsilon^{2}.
\end{align}
Thus, we have $\|\partial_{\omega_{j}}u\|_{\sigma,s}\le K_{4}\varepsilon^{2}$ by choosing appropriate $K_{4}$.
Similarly, we have $\|\partial_{\bm{a}}u\|_{\sigma,s}\le K_{4}\varepsilon^{2}$. 
\end{proof}

\section{The bifurcation equations}
Substituting $u(\omega,v)$ into the bifurcation equation \eqref{bifurcation}, we obtain
\begin{equation}
    QT_{\omega}Qu+Q\Delta(\varepsilon^{-1}S(\varepsilon v+\varepsilon u(\omega,v)))=0.
\end{equation}
Passing to Fourier coefficients, we obtain
\begin{equation}
    (-k\cdot\omega+|n_{k}|^{2}+\rho)a_{k}+\varepsilon^{-1}|n_{k}|^{2} (S(\varepsilon v+\varepsilon u(\omega,v)))^{\land}(k,n_{k})=0,\ \ (k,n_{k})\in\mathcal{R}.
\end{equation}
Note that $S(u)=2|u|^{2}u+O(|u|^{4})$ and $u(\omega,v)=O(\varepsilon^{2})$.
We have
\begin{equation}
    (-k\cdot\omega+|n_{k}|^{2}+\rho)a_{k}+2\varepsilon^{2}|n_{k}|^{2} (|v|^{2}v)^{\land}(k,n_{k})+O(\varepsilon^{3})=0,\ \ (k,n_{k})\in\mathcal{R}.
    \label{bifurcationzhenggebianhuanhou}
\end{equation}

\subsection{The structure of the resonant set $\mathcal{R}$}\label{structureofR}

Recall 
\begin{equation}
    \mathcal{R}=\{(k,n)|-k\cdot\omega_{0}+|n|^{2}+\rho=0\}\cap \mathcal{W}.
\end{equation}
Before solving the bifurcation equation, we first explore the structure of $\mathcal{R}$.

Obviously, $(e_{j},n_{j})\in\mathcal{R}$. By notation abuse, we denote $k\in\mathcal{R}$ if $(k,n_{k})\in\mathcal{R}$. When $|k|\gg 1$, we have $|-k\cdot \omega+|n_{k}|^{2}+\rho|>0$. Thus, $|\mathcal{R}|<\infty$. Denote $\textup{supp }k=\{j|k_{j}\ne0\}$.

\begin{lemma}\label{yinli4.1}
    Assume that $\rho$ is irrational and the vectors $n_1, \dots, n_b$ are linearly independent. Then, for any $k \in \mathcal{R}$, either
\begin{equation}
k = e_j \quad \text{for some } j,
\end{equation}
or
\begin{equation}
|\operatorname{supp} k| \ge 3.
\label{kdetiaojian1}
\end{equation}

If, in addition, for any distinct indices $i, j, l$, the vectors $n_i - n_j$ and $n_i - n_l$ are not perpendicular, i.e., $(n_i - n_j) \not\perp (n_i - n_l)$, then condition \eqref{kdetiaojian1} is strengthened to
\begin{equation}
|\operatorname{supp} k| \ge 3 \quad \text{but} \quad k \ne \pm e_i \pm e_j \pm e_l \quad \text{for any distinct } i, j, l.
\label{kdetiaojian2}
\end{equation}

Furthermore, if   for any distinct $i, j, l$, the angle between $n_i - n_j$ and $n_i - n_l$ is less than $\frac{\pi}{2}$, then \eqref{kdetiaojian1} is strengthened to
\begin{equation}
|\operatorname{supp} k| \ge 4.
\end{equation}
\end{lemma}
\begin{proof}
We have
\begin{equation}
    -k\cdot\omega_{0}+|n_{k}|^{2}+\rho=-\sum_{j=1}^{b}k_{j}|n_{j}|^{2}+|\sum_{j=1}^{b}k_{j}n_{j}|^{2}+(1-\sum_{j=1}^{b}k_{j})\rho.
\end{equation}
    Since $\rho$ is irrational, $(k,n_{k})\in\mathcal{R}$ if and only if 
    \begin{align}
        &\sum_{j=1}^{b}k_{j}=1,\\
        &\sum_{j=1}^{b}k_{j}|n_{j}|^{2}=|\sum_{j=1}^{b}k_{j}n_{j}|^{2},
    \end{align}
    which is equivalent to
    \begin{align}
         &\sum_{j=1}^{b}k_{j}=1,\label{chaopingmian}\\
        &\sum_{j\ne l}k_{j}k_{l}|n_{j}-n_{l}|^{2}=0.\label{zhijipanduan}
    \end{align}
    This comes from the equality 
    \begin{equation}
        \sum_{j=1}^{b}\alpha_{j}|x_{j}|^{2}-|\sum_{j=1}^{b}\alpha_{j}x_{j}|^{2}=\frac{1}{2}\sum_{j\ne l}\alpha_{j}\alpha_{l}|x_{j}-x_{l}|^{2},
        \end{equation}
        when $\sum_{j=1}^{b}\alpha_{j}=1$.
    If $|\textup{supp }k|=1$, we obtain $k=e_{j}$ for some $j$ from \eqref{chaopingmian}. Moreover, it is straightforward to deduce a contradiction between  $|\textup{supp }k|=2$ and \eqref{zhijipanduan}. 

    Now, assume that for any distinct indices $i, j, l$, the vectors $n_i - n_j$ and $n_i - n_l$ are not perpendicular.
    Suppose $k = \pm e_i \pm e_j \pm e_l  \text{ for some distinct } i, j, l$. Without loss of generality, suppose $|\textup{supp }k|=3$.
    We have
    \begin{equation}
        |k_{1}k_{2}|n_{1}-n_{2}|^{2}+k_{2}k_{3}|n_{2}-n_{3}|^{2}+k_{1}k_{3}|n_{1}-n_{3}|^{2}|=\nu_{1}|n_{1}-n_{2}|^{2}+\nu_{2}|n_{2}-n_{3}|^{2}+\nu_{3}|n_{1}-n_{3}|^{2},
    \end{equation}
    where $\#\{\nu_{j}|\nu_{j}=1\}=2$. This is a contradiction with the assumption.
    
    Now, assume that for any distinct $i, j, l$, the angle between $n_i - n_j$ and $n_i - n_l$ is less than $\frac{\pi}{2}$. Suppose $|\textup{supp }k|=3$. Without loss of generality, assume $\textup{supp }k=\{k_{1},k_{2},k_{3}\}$. Since $k_{1}+k_{2}+k_{3}=1$, there is at least a $k_{j}>0$ and a $k_{j'}<0$. Suppose $k_{1}>0$ and $k_{2}<0$. 

    \textbf{Case 1}. $k_{3}>0$.
    
    We have
    \begin{equation}
        |k_{2}|=k_{1}+k_{3}-1\ge\max\{k_{1},k_{3}\}.
    \end{equation}
    Thus, we have
    \begin{align}
        &|k_{1}k_{2}|n_{1}-n_{2}|^{2}+k_{2}k_{3}|n_{2}-n_{3}|^{2}+k_{1}k_{3}|n_{1}-n_{3}|^{2}|\notag\\
        \ge & |k_{1}k_{2}|\cdot|n_{1}-n_{2}|^{2}+|k_{2}k_{3}|\cdot|n_{2}-n_{3}|^{2}-|k_{1}k_{3}|\cdot |n_{1}-n_{3}|^{2}\notag\\
        \ge & |k_{1}k_{3}|\cdot(|n_{1}-n_{2}|^{2}+|n_{2}-n_{3}|^{2}-|n_{1}-n_{3}|^{2})>0,
    \end{align}
    which contradicts with \eqref{zhijipanduan}.

     \textbf{Case 2}. $k_{3}<0$.

     We have
     \begin{equation}
     k_{1}=1+|k_{1}|+|k_{2}|.    
     \end{equation}
     Thus, we have
     \begin{align}
        &|k_{1}k_{2}|n_{1}-n_{2}|^{2}+k_{2}k_{3}|n_{2}-n_{3}|^{2}+k_{1}k_{3}|n_{1}-n_{3}|^{2}|\notag\\
        \ge & |k_{1}k_{2}|\cdot|n_{1}-n_{2}|^{2}+|k_{2}k_{3}|\cdot|n_{2}-n_{3}|^{2}-|k_{1}k_{3}|\cdot |n_{1}-n_{3}|^{2}\notag\\
        \ge & |k_{1}k_{3}|\cdot(|n_{1}-n_{2}|^{2}+|n_{2}-n_{3}|^{2}-|n_{1}-n_{3}|^{2})>0,
    \end{align}
    which contradicts  \eqref{zhijipanduan}.
\end{proof}

We give another collection of $n_{j}$ such that $\mathcal{R}$ only contains $\{e_{j}\}$.
\begin{lemma}\label{yinli4.2}
    Suppose that $|n_{j}|^{2}\equiv M$ for some constant $M>0$, and that $n_{j}\perp n_{l}$ for all $j\ne l$. Then, $\mathcal{R}=\{(e_{j},n_{j})|j=1,...,b\}$.
\end{lemma}
\begin{proof}
    It suffices to prove that the solution set of

\begin{equation}
\begin{cases}
k \in \mathbb{Z}^{b}, \\
\displaystyle \sum_{j=1}^{b} k_{j} = 1, \\
\displaystyle \sum_{j=1}^{b} k_{j} |n_{j}|^{2} = \left| \sum_{j=1}^{b} k_{j} n_{j} \right|^{2}
\end{cases}\notag
\end{equation}
    is $\{e_{j}|1\le j\le b\}$. In fact, we have
    \begin{equation}
        \left|\sum_{j=1}^{b}k_{j}n_{j}\right|^{2}=M\sum_{j=1}^{b}k_{j}^{2}.
    \end{equation}
    On the other hand, we have
    \begin{equation}
        \sum_{j=1}^{b}k_{j}|n_{j}|^{2}=M\sum_{j=1}^{b}k_{j}.
    \end{equation}
    Thus, we have
    \begin{align}
        \sum_{j=1}^{b}k_{j}^{2}=\sum_{j=1}^{b}k_{j}=1,
    \end{align}
    which implies $k=e_{j}$, $j=1,...,b$.
\end{proof}

\subsection{Solving the bifurcation equations}
In this section, we suppose that $k\ne \pm e_{i}\pm e_{j}\pm e_{l}$ for any $k\in\mathcal{R}$.

We divide the bifurcation equations \eqref{bifurcationzhenggebianhuanhou} into two parts:
\begin{align}
    &(-\omega_{j}+|n_{j}|^{2}+\rho)a_{j}+2\varepsilon^{2}|n_{j}|^{2}(|v|^{2}v)^{\land}(e_{j},n_{j})+O(\varepsilon^{3})=0,\ \ j=1,...,b,\label{Q1}\\
    &(-k\cdot\omega+|n_{k}|^{2}+\rho)a_{k}+2\varepsilon^{2}|n_{k}|^{2} (|v|^{2}v)^{\land}(k,n_{k})+O(\varepsilon^{3})=0,\ \ k\notin\{e_{1},...,e_{b}\}.\label{Q2}
\end{align}
We call \eqref{Q1} $Q$-I equations and call \eqref{Q2} $Q$-II equations.
From now on, we denote by $a_{k}$ the $a_{k}$ where $k\notin\{e_{1},...,e_{b}\}$. Moreover, we denote $\bm{a}_{k}=(a_{k}:k\ne e_{j}, 1\le j\le b)$, $\bm{a}=(a_{1},...,a_{b},\bm{a}_{k})$ and $\tilde{b}=|\mathcal{R}|-b$. 

\subsubsection{Q-I equations}
First, we consider $Q$-I equations. We have the following statement.

\begin{proposition}\label{mingti4.3}
Let $a_{k}\in[-\varepsilon,\varepsilon]$ and $a_{j}\in[1,2]$. There exists a solution $\omega(\bm{a})$ to \eqref{Q1} such that
\begin{equation}
    \omega_{j}=|n_{j}|^{2}+\rho+\varepsilon^{2}|n_{j}|^{2}\left(4\sum_{l=1}^{b}a_{l}^{2}-2a_{j}^{2}\right)+h_{j}(\bm{a}),
\end{equation}
where $h_{j}=O(\varepsilon^{3})$ and $\frac{\partial h_{j}}{\partial \bm{a}}=O(\varepsilon^{3})$.
\end{proposition}
\begin{proof}
    We aim to solve $\omega(\bm{a})$ from $Q$-I equations. Note that we have
\begin{align}
    (|v|^{2}v)^{\land}(e_{j},n_{j})&=\int v^{2}\bar{v}e^{-i(e_{j},n_{j})\cdot(\theta,x)}d\theta d x\notag\\
    &=\sum_{\substack{k+k'-k''=e_{j}}} a_{k}a_{k'}a_{k''}\notag\\
    &=2\sum_{l\ne j}a_{l}^{2}a_{j}+a_{j}^{3}+O(|a_{k}|^{2}).
\end{align}
Thus, the $Q$-I equations \eqref{Q1} become
\begin{equation}
    (-\omega_{j}+|n_{j}|^{2}+\rho)a_{j}+2\varepsilon^{2}|n_{j}|^{2}\left(2\sum_{l\ne j}a_{l}^{2}a_{j}+a_{j}^{3}+O(|a_{k}|^{2})\right)+O(\varepsilon^{3})=0,\ \ j=1,...,b.\label{dairuQ1}
\end{equation}
Denote 
\begin{equation}
    \sigma_{j}=\frac{-\omega_{j}+|n_{j}|^{2}+\rho}{\varepsilon^{2}}.
    \label{sigmaomegabianhuan}
\end{equation}
Then, \eqref{dairuQ1} becomes
\begin{equation}
    \sigma_{j}a_{j}+2|n_{j}|^{2}\left(\left(2\sum_{l=1}^{b}a_{l}^{2}-a_{j}^{2}\right)a_{j}+O(|a_{k}|^{2})\right)+O(\varepsilon)=0,\ \ j=1,...,b.\label{Q1zuizhongban}
\end{equation}
Moreover, we have
\begin{equation}
    \sigma_{j}+2|n_{j}|^{2}\left(2\sum_{l=1}^{b}a_{l}^{2}-a_{j}^{2}\right)+O(|a_{k}|^{2})+O(\varepsilon)=0,\ \ j=1,...,b.\label{Q1zuizhongban1}
\end{equation}
We refer to the equation corresponding to a given index $j$ as the $j$th-equation.
Denote the left side of \eqref{Q1zuizhongban1} by $g_{j}(\sigma(\omega),\bm{a})$ and the $O(\varepsilon)$ term by $\tilde{g}_{j}(\omega,\bm{a})$. Furthermore, we denote
\begin{equation}
    f_{j}(a_{1},...,a_{b})=2|n_{j}|^{2}\left(2\sum_{l=1}^{b}a_{l}^{2}-a_{j}^{2}\right).
\end{equation}
By \eqref{sigmaomegabianhuan}, the definition of $\tilde{g}_{j}$ and (iv) of Lemma \ref{rangesol}, we have
\begin{equation}
    \frac{\partial \tilde{g}}{\partial\sigma}=\frac{\partial \tilde{g}}{\partial\omega}\cdot\frac{\partial\omega}{\partial\sigma}=O(\varepsilon^{3}),\ \ \ \ \  \frac{\partial\tilde{g}}{\partial \bm{a}}=O(\varepsilon).
\end{equation}
We solve \eqref{Q1zuizhongban1} from $j=1$ to $j=b$. Regard $\sigma_{2},...,\sigma_{b}\in [-\frac{1}{\varepsilon},\frac{1}{\varepsilon}]$ and $\bm{a}$ as parameter. We have the solution to $1$th-equation as follows
\begin{equation}
    \sigma_{1}(\sigma_{2},...,\sigma_{b},\bm{a})=-f_{1}(a_{1},...,a_{b})+\tilde{h}_{1}(\sigma_{2},...,\sigma_{b},\bm{a}),
\end{equation}
where $\tilde{h}_{1}=O(\varepsilon)$. By the implicit function theorem, we have
\begin{align}
    &\frac{\partial \sigma_{1}}{\partial \sigma_{j}}=-\left(\frac{\partial g_{1}}{\partial \sigma_{1}}\right)^{-1}\frac{\partial g_{1}}{\partial \sigma_{j}}=O(\varepsilon^{3}), \ \ \textup{ for } j\ne 1,\\
    &\frac{\partial\sigma_{1}}{\partial \bm{a}_{k}}=-\left(\frac{\partial g_{1}}{\partial \sigma_{1}}\right)^{-1}\frac{\partial g_{1}}{\partial \bm{a}_{k}}=O(\varepsilon).
\end{align}
Furthermore, we have
\begin{align}
    \frac{\partial\sigma_{1}}{\partial a_{j}}
    &=-\left(\frac{\partial g_{1}}{\partial a_{j}}\right)^{-1}\frac{\partial g_{1}}{\partial a_{j}}\notag\\
    &=-(1+O(\varepsilon^{3}))^{-1}\left(\frac{\partial f_{1}}{\partial a_{j}}+O(\varepsilon)\right)\notag\\
    &=-\frac{\partial f_{1}}{\partial a_{j}}+O(\varepsilon), \textup{ for }j=1,...,b,\label{ajderiv}\\
\end{align}
which implies that
\begin{equation}
    \frac{\partial \tilde{h}_{1}}{\partial a_{j}}=O(\varepsilon),\ \  \textup{ for } j=1,...,b.
\end{equation}
Substituting $\sigma_{1}(\sigma_{2},...,\sigma_{b},\bm{a})$ into the $2$th-equation and following the same argument, we obtain
\begin{equation}
    \sigma_{2}(\sigma_{3},...,\sigma_{b},\bm{a})=-f_{2}(a_{1},...,a_{b})+\tilde{h}_{2}(\sigma_{3},...,\sigma_{b},\bm{a}),
\end{equation}
where
\begin{align}
    &\tilde{h}_{2}=O(\varepsilon),\\
    &\frac{\partial\tilde{h}_{2}}{\partial\sigma_{j}}=O(\varepsilon^{3}),\ \ \textup{ for } j\ge 3,\\
    &\frac{\partial \tilde{h}_{2}}{\partial\bm{a}}=O(\varepsilon).
\end{align}
Iterate this process till $j=b$. We obtain the solution $\bm{\sigma}(\bm{a})$ of \eqref{Q1zuizhongban1} satisfying
\begin{align}
    &\sigma_{j_{0}}=-f_{j_{0}}(a_{1},...,a_{b})+\tilde{h}_{j_{0}}(\sigma_{j_{0}+1},...,\sigma_{b},\bm{a}),\\
    &\tilde{h}_{j_{0}}=O(\varepsilon),\\
    &\frac{\partial \tilde{h}_{j_{0}}}{\partial\sigma_{j}}=O(\varepsilon^{3}), \textup{ for }j\ge j_{0}+1,\\
    &\frac{\partial \tilde{h}_{j_{0}}}{\partial\bm{a}}=O(\varepsilon).
\end{align}
Successive substitution of $\sigma_{j+1}$ into $\sigma_{j}$, combined with \eqref{sigmaomegabianhuan}, yields the desired result.
\end{proof}
Note that if there in no $Q$-II equations, we have solved the bifurcation equation. Thus, we obtain Theorem \ref{dingli1.1}. Moreover, combing  Theorem \ref{dingli1.1}, Lemma \ref{yinli4.1} and Lemma \ref{yinli4.2}, we obtain Corollary \ref{tuilun1.2}.

\subsubsection{Q-II equations}
In this section, we consider the $Q$-II equations.

Let $k_{0}\in\mathcal{R}\setminus\{(e_{j}, n_{j})\}$. For simplicity, we denote $\{e_{j}:1\le j\le b\}$ by $\{e_{j}\}$.

We have
\begin{align}
    (|v|^{2}v)^{\land}(e_{j},n_{j})&=\int v^{2}\bar{v}e^{-i(e_{j},n_{j})\cdot(\theta,x)}d\theta d x\notag\\
    &=\sum_{\substack{k+k'-k''=k_{0}}} a_{k}a_{k'}a_{k''}\notag\\
    &=\left(2\sum_{j=1}^{b}a_{j}^{2}\right)a_{k_{0}}+\sum_{k\ne k_{0}}A_{k_{0} k}(a_{1},...,a_{b})a_{k}+O(|\bm{a}_{\tilde{k}}|^{2};\tilde{k}\notin\{e_{j}\}),
\end{align}
where $A_{k_{0}k}(a_{1},...,a_{b})$ is a homogeneous quadratic polynomial in $a_{j}$ and  only depends on $n_{j}$. 
Note that we have
\begin{align}
    &-k\cdot\omega+|n_{k}|^{2}+\rho\notag\\
    =&-\sum_{j=1}^{b}k_{j}\left(|n_{j}|^{2}+\rho+\varepsilon^{2}|n_{j}|^{2}\left(4\sum_{l=1}^{b}a_{l}^{2}-2a_{j}^{2}\right)+O(\varepsilon^{3})\right)+|n_{k}|^{2}+\rho\notag\\
    =&-\varepsilon^{2}\sum_{j=1}^{b}k_{j}|n_{j}|^{2}\left(4\sum_{l=1}^{b}a_{l}^{2}-2a_{j}^{2}\right)\notag\\
    =&-\varepsilon^{2}|n_{k}|^{2}\cdot4\sum_{l=1}^{b}a_{l}^{2}+\varepsilon^{2}\sum_{j=1}^{b}2k_{j}a_{j}^{2}|n_{j}|^{2}.
\end{align}

Substituting $\omega(\bm{a})$ into the $Q$-II equations, we obtain
\begin{equation}
    \left( \sum_{j=1}^{b}k_{j}a_{j}^{2}|n_{j}|^{2} \right)a_{k}+\sum_{k'\ne k}A_{kk'}a_{k'}+O(|\bm{a}_{\tilde{k}}|^{2};\tilde{k}\notin\{e_{j}\})+O(\varepsilon)=0.
    \label{zhongjianbanbenQ2}
\end{equation}
Denote $A_{kk}=\sum_{j=1}^{b}k_{j}a_{j}^{2}|n_{j}|^{2}$. Note that $A_{kk}$ also only depends on $n_{j}$. 

Thus $A=(A_{kk'})_{k,k'\in\mathcal{R}\setminus\{(e_{j},n_{j})\}}$  is a matrix determined by $n_{j}$, whose entries are homogeneous quadratic 
in $a_{j}$. Moreover, the coefficients of $A_{kk'}(a_{1},...,a_{b})$ are integers. With the notation $A$, \eqref{zhongjianbanbenQ2} can be written as 
\begin{equation}
    A\bm{a}_{k}+O(|\bm{a}_{k}|^{2})+O(\varepsilon)=0.
    \label{zuizhongbanQ2}
\end{equation}
\begin{proposition}\label{mingti4.4}
    Suppose that $\det A \not\equiv 0$. There exists a subset $I_{\varepsilon}\subset[1,2]^{b}$ such that
    \begin{itemize}
        \item[\textup{(i)}] $\textup{mes } I_{\varepsilon}<\varepsilon^{c}$, where $c$ is a constant depending on $b$, $\tilde{b}$.
        \item[\textup{(ii)}] For all $(a_{1},...,a_{b})\in [1,2]^{b}\setminus I_{\varepsilon}$, there exists a solution $\bm{a}_{k}(a_{1},...,a_{b})$ to \eqref{zuizhongbanQ2}. Moreover, we have $\bm{a}_{k}(a_{1},...,a_{b})=O(\varepsilon^{\frac{3}{4}})$.
    \end{itemize}
\end{proposition}
\begin{proof}
Define
\begin{equation}
    I_{\varepsilon}=\{(a_{1},...,a_{b}):|\det A|\le \varepsilon^{\frac{1}{6}}\}.
\end{equation}
Note that  $\det A\in\mathbb{Z}[a_{1},...,a_{b}]$, $\det A\not\equiv 0$ and the degree of $\det A$ is less than or equal to $2\tilde{b}$. By Lemma A.4 in \cite{xyy2025construction}, we have
\begin{equation}
    \textup{mes } I_{\varepsilon}<C_{\tilde{b}}b\varepsilon^{\frac{1}{6}\cdot\frac{1}{2b\tilde{b}}}.
\end{equation}
Let $(a_{1},...,a_{b})\in[1,2]^{b}\setminus I_{\varepsilon}$. By Cramer's rule, we have
\begin{equation}
    |(A^{-1})_{kk'}|<\frac{C(n_{1},...,n_{b})}{|\det A|}<C(n_{1},...,n_{b})\varepsilon^{-\frac{1}{6}}.
\end{equation}
    Thus, $A^{-1}$ exists and 
\begin{equation}
    \|A^{-1}\|<\varepsilon^{-\frac{1}{5}}
\end{equation}
for  $(a_{1},...,a_{b})\in[1,2]^{b}\setminus I_{\varepsilon}$. Denote the last two terms of \eqref{zuizhongbanQ2} by $r(a_{1},...,a_{b},\bm{a}_{k})$. Then, the equation \eqref{zuizhongbanQ2} is equivalent to
\begin{equation}
    \bm{a}_{k}=-A^{-1}r(a_{1},...,a_{b},\bm{a}_{k}).
\end{equation}
Denote $\mathcal{G}'(\bm{a}_{k})=-A^{-1}r(a_{1},...,a_{b},\bm{a}_{k})$. We now prove $\mathcal{G}'$ is a contraction map on $B(0,\varepsilon^{\frac{3}{4}})$.

Suppose that $\bm{a}_{k}\in B(0,\varepsilon^{\frac{3}{4}})$. We have
\begin{equation}
    r(a_{1},...,a_{b},\bm{a}_{k})=O(\varepsilon^{\frac{3}{2}})+O(\varepsilon)=O(\varepsilon),
\end{equation}
which implies
\begin{align}
    \|\mathcal{G}'(\bm{a}_{k})\|\le\|A^{-1}\|\cdot\|r(a_{1},...,a_{b},\bm{a}_{k})\|<\varepsilon^{\frac{3}{4}}.
\end{align}
On the other hand, suppose $\bm{a}_{k}', \bm{a}_{k}''\in B(0,\varepsilon^{\frac{3}{4}})$, we have
\begin{align}
    \|\mathcal{G}'(\bm{a}_{k}')-\mathcal{G}'(\bm{a}_{k}'')\|&=\left\|A^{-1}\int_{0}^{1}\frac{\partial r(a_{1},...,a_{b},\bm{a}_{k}+t(\bm{a}_{k}''-\bm{a}_{k}'))}{\partial \bm{a}_{k}}\cdot(\bm{a}_{k}''-\bm{a}_{k}')dt \right\|\notag\\
    &\le C\varepsilon^{-\frac{1}{5}}\varepsilon\|\bm{a}_{k}'-\bm{a}_{k}''\|\notag\\
    &<\frac{1}{2}\|\bm{a}_{k}'-\bm{a}_{k}''\|.
\end{align}
Thus, there exists an unique solution $\bm{a}_{k}(a_{1},...,a_{b})=O(\varepsilon^{\frac{3}{4}})$ to \eqref{zuizhongbanQ2} (i.e., the $Q$-II equation).
\end{proof}
Combing Lemma \ref{rangesol}, Lemma \ref{yinli4.1},  Proposition \ref{mingti4.3} and Proposition \ref{mingti4.4}, we obtain Theorem \ref{dingli1.3}.
\section{Proof of Theorem \ref{dingli1.4}}
Suppose that $H(u,\bar{u})$ takes the form
\begin{equation}
    H(u,\bar{u})=|u|^{4}+\sum_{p=3}^{\infty}\alpha_{p}|u|^{2p}.\notag
\end{equation}
Replacing $u=e^{-i\rho' t}\tilde{u}$ into \eqref{yuanshifangcheng}, we obtain
\begin{equation}
    i\tilde{u}_{t}-\Delta \tilde{u}+(\rho+\rho')\tilde{u}+\Delta (\frac{\partial H}{\partial \bar{u}}(\tilde{u},\bar{\tilde{u}}))=0.
\end{equation}
Choose $\rho'$ such that $\rho+\rho'$ satisfies the Diophantine conditions. Applying Theorem \ref{dingli1.1}, Corollary \ref{tuilun1.2} and Theorem \ref{dingli1.3}, we obtain Theorem \ref{dingli1.4}.

\
\

\noindent \textbf{Data Availability} No datasets were generated or analyzed during the current study.

\
\

\noindent \textbf{Conflict  of Interest} The authors have no conflict of interest.

\bibliographystyle{plainurl}
\bibliography{reference}

\end{document}